\documentclass[12pt]{article}
\usepackage{amssymb}
\usepackage{amsmath}
\usepackage{multicol}
\newcommand{\lap}{\mbox{$\bigtriangleup$}}

\newcommand{\be}{\begin{equation}}
\newcommand{\ee}{\end{equation}}

\newtheorem{mthm}{Theorem}

\newtheorem{mlem}{Lemma}

\newtheorem{pro}{Proposition}[section]
\newtheorem{lem}{Lemma}[section]

\begin{document}

\title{ The A Priori Estimate and Existence of the Positive Solution for A Nonlinear System Involving the Fractional  Laplacian }
\author {Ran Zhuo \thanks{Partially supported by NSFC 11701207.}  \hspace{.2in} Yan Li\thanks{Corresponding author. }}

\date{}
\maketitle

\begin{abstract}

In the paper, we consider the fractional
elliptic system
\begin{equation*}\left\{\begin{array}{ll}
(- \Delta)^{\frac{\alpha_1}{2}}u(x)+\sum\limits^n_{i=1}b_i(x)\frac{\partial u}{\partial x_i}+B(x)u(x)=f(x,u,v),& \mbox {in\,\,\,} \Omega,\\
(- \Delta)^{\frac{\alpha_2}{2}}v(x)+\sum\limits^n_{i=1}c_i(x)\frac{\partial
v}{\partial x_i}+C(x)v(x)=g(x,u,v),& \mbox {in\,\,\,}
\Omega,\\
u=v=0, & \mbox {in\,\,\,} \mathbb{R}^n\setminus\Omega,
\end{array} \right.\label{a-1.2}
\end{equation*}
where $\Omega$ is a bounded domain with $C^2$ boundary in $\mathbb{R}^n$ and $n>\max\{\alpha_1,\alpha_2\}$. We first utilize the blowing-up and re-scaling method to derive the a priori estimate for positive solutions when
$1<\alpha_1,\alpha_2 <2$. Then for $0<\alpha_1,\alpha_2
<1$, we obtain the regularity estimate of positive solutions. On top of this, using the topological degree theory we prove the existence of positive solutions.

\end{abstract}

\maketitle

\section{Introduction}

The fractional Laplacian in $\mathbb{R}^n$ is a nonlocal pseudo-differential
operator defined as
\begin{equation}
(-\Delta)^{\alpha/2} u(x) := C_{n,\alpha} \, PV \int_{\mathbb{R}^n}
\frac{u(x)-u(z)}{|x-z|^{n+\alpha}} dz, \label{Ad7}
\end{equation}
 where $\alpha$ is any real number between $0$ and $2$ and PV stands for the Cauchy principal value. Let
$$L_{\alpha}:=\{u: \mathbb{R}^n\rightarrow \mathbb{R} \mid \int_{\mathbb{R}^n}\frac{|u(x)|}{1+|x|^{n+\alpha}} \, d x <\infty\}.$$
Then for $u \in L_{\alpha}\cap C^{1,1}_{loc}(\Omega)$, one can see that the operator in (\ref{Ad7}) is well-defined through elementary calculation.  For readers who are interested in more discussions on the definitions of the fractional Laplacian, we refer them to  \cite{Si} and \cite{CLM}.

The fractional Laplacian operator $(-\Delta)^{\alpha/2}$  appears in a wide class of physical systems, including L\'{e}vy flights and stochastic interfaces(\cite{ZRK}).
Due to its non-locality as can be seen in (\ref{Ad7}), the fractional Laplacian and many more general nonlocal operators have been studied in multiple areas, such as optimization(\cite{CRS}), flame propagation(\cite{DL}),  finance(\cite{CT})
and  phase transitions(\cite{ABS}, \cite{G}). The fractional Laplacian  can also be
seen as the infinitesimal generator of a L\'{e}vy process and has been applied in probability(\cite{B}). As a result of its wide applications, fruitful results have been obtained in the properties of solutions to systems involving the fractional Laplacian.

In \cite{QX}, for $0<\alpha<1$ the authors derived nonexistence of positive bounded viscosity solution for
$$
\left\{ \begin{array}{ll}
(-\triangle)^{\alpha} u(x) =v^p(x), & x\in\mathbb{R}^n,  \\
(-\triangle)^{\alpha }v (x) = u^q(x) ,&  x\in \mathbb{R}^n,
\end{array} \right.
$$
when $p,\, q$ take values in certain range.

For $0<\alpha, \, \beta<2$,  symmetry of solutions was derived in \cite{LM}  for
$$
\left\{ \begin{array}{ll}
(-\triangle)^{\alpha/2} u(x) =f(v(x)), & x\in\mathbb{R}^n,  \\
(-\triangle)^{\beta/2 }v (x) =g( u(x)) ,&  x\in \mathbb{R}^n,  \\
u, \, v \geq 0,  ,&  x\in \mathbb{R}^n,
\end{array} \right.
$$
under some bounded-ness assumptions about the nonlinear terms.

In \cite{DP}, for $s_1\,, s_2 \in (0,1)$ the authors studied
$$
\left\{ \begin{array}{ll}
(-\triangle)^{s_1} u(x) =F_1(u, v), & x\in\mathbb{R}^n,  \\
(-\triangle)^{s_2 }v (x) = F_2(u, v) ,&  x\in \mathbb{R}^n.
\end{array} \right.
$$
Under the assumption that $F_i \in C_{loc} ^{1,1}(\mathbb{R}^2)$ they obtained some Poincar\'{e}-type formulas for the $\alpha$-harmonic extension in the half-space, with which they continued to prove a symmetry result both for stable and for monotone solutions

For $0<\alpha<2$, the authors in \cite{ZCCY} proved the non-existence in the subcritical case for positive solutions of
$$
\left\{ \begin{array}{ll}
(-\triangle)^{\alpha/2} u_i(x) =f_i(u_1(x),  \cdots, u_m(x)), & x\in\mathbb{R}^n,  \\
u_i(x)\geq 0  ,&  x\in\mathbb{R}^n,
\end{array} \right.
$$
by studying its equivalent integral system.

Recently,the author in \cite{L} considered
$$
\left\{\begin{array}{ll}
(-\triangle)^{\frac{\alpha}{2}}u+\sum^{N}_{i = 1}b_{i}(x)\frac{\partial u}{\partial x_{i}}+C(x)u = f(x,v), & x\in \Omega,\\
(-\triangle)^{\frac{\beta}{2}}v+\sum^{N}_{i = 1}c_{i}(x)\frac{\partial v}{\partial x_{i}}+D(x)v = g(x,u),&  x\in \Omega,\\
u>0,\, v>0, & x\in \Omega,\\
u = 0, v = 0, & x\in \mathbb R^{N}\setminus \Omega,
\end{array}\right.
$$
for  $\alpha, \beta \in (1,2)$. Using the scaling method, the author obtained the a priori bounds of positive solutions by imposing certain conditions on potential functions and the nonlinear terms.

Motivated by the works above, we are interested in studying the following fractional system, which can be seen as a generalized version of the equations listed above. We consider
\begin{equation}\left\{\begin{array}{ll}
(- \Delta)^{\frac{\alpha_1}{2}}u(x)+\sum\limits^n_{i=1}b_i(x)\frac{\partial u}{\partial x_i}+B(x)u(x)=f(x,u,v),& \mbox {in\,\,\,} \Omega,\\
(- \Delta)^{\frac{\alpha_2}{2}}v(x)+\sum\limits^n_{i=1}c_i(x)\frac{\partial
v}{\partial x_i}+C(x)v(x)=g(x,u,v),& \mbox {in\,\,\,}
\Omega,\\
u=v=0, & \mbox {in\,\,\,} \mathbb{R}^n\setminus\Omega.
\end{array} \right.\label{In-1.2}
 \end{equation}
In order to make the first-order derivative terms meaningful, we
introduce a natural restriction  $1<\alpha_1,\, \alpha_2<2$. In this paper
we assume $\Omega$ to be a  bounded  domain in $\mathbb{R}^n$ with  $C^2$ boundary and
$n>\max\{\alpha_1,\, \alpha_2\}$.

Below is the collection of assumptions we will impose on the coefficients in (\ref{In-1.2}) for  the discussion of regularity.
\begin{description}
  \item[$(A_1)$]
  $b_i(x)$, $c_i(x)$, $B(x)$ and $C(x)$
are  bounded on $\bar{\Omega}$.
  \item[$(A_2)$]
  $f(x,s,t),\, g(x,s,t) \in C(\bar{\Omega}, R, R)$.
  \item[$(A_3)$]
  For positive $p_i,\,q_i$, $i=1,\, 2$, satisfying $p_1+\frac{n-\alpha_2}{n-\alpha_1}q_1<\frac{n+\alpha_1}{n-\alpha_1}$, \,
$ p_2+\frac{ n-\alpha_2}{n-\alpha_1 } q_2< \frac{n+\alpha_2}{n-\alpha_1 }$,
assume $$\lim\limits_{\max\{|s|,|t|\}\to
\infty}\frac{f(x,s,t)}{|s|^{p_1}|t|^{q_1}}=k_1(x), \quad
\lim\limits_{\max\{|s|,|t|\}\to
\infty}\frac{g(x,s,t)}{|s|^{p_2}|t|^{q_2}}=h_1(x),$$
where
$k_1(x), \,h_1(x)$ are positive functions continuous on $
\bar{\Omega}$.
  \item[$(A_4)$]
  For positive $\tau_i,\, \eta_i$, $i=1,\, 2$, assume $\tau_i+\eta_i>1$, $i=1,\, 2$,
 $$\lim\limits_{\min\{|s|,|t|\}\to
0}\frac{f(x,s,t)}{|s|^{\tau_1}|t|^{\eta_1}}=k_2(x), \quad
\lim\limits_{\min\{|s|,|t|\}\to
0}\frac{g(x,s,t)}{|s|^{\tau_2}|t|^{\eta_2}}=h_2(x),$$
 where $k_2(x), \,h_2(x)$ are positive functions continuous on $
\bar{\Omega}$.
\end{description}

 Let
$d(x):=dist(x,\partial\Omega)$, for  $\theta>0$, we define
$$\|u\|_{0,\theta}:=\sup_{\Omega}\big(d^\theta(x)|u(x)|\big),$$
$$\|u\|_{1,\theta}:=\sup_{\Omega}\big(d^\theta(x)|u(x)|+d^{\theta+1}(x)\nabla u(x)\big).$$

We first give the following regularity estimate for (\ref{In-1.2}).

\begin{mlem}
For $1<\alpha_1,\,\alpha_2<2$, suppose $(A_1)-(A_3)$ hold. If $u,\,v\in C^1(\Omega)\cap L^{\infty}(R^n)$ satisfy the following system
\begin{equation}\left\{\begin{array}{ll}
(- \Delta)^{\frac{\alpha_1}{2}}u(x)+\sum\limits^n_{i=1}b_i(x)\frac{\partial u}{\partial x_i}+B(x)u(x)=f(x,u,v),& \mbox {in\,\,\,} \Omega,\\
(- \Delta)^{\frac{\alpha_2}{2}}v(x)+\sum\limits^n_{i=1}c_i(x)\frac{\partial
v}{\partial x_i}+C(x)v(x)=g(x,u,v),& \mbox {in\,\,\,}
\Omega,\\
u(x), \, v(x)>0, & \mbox {in\,\,\,} \Omega.
\end{array} \right.\label{lem-2-1} \end{equation}
Then there exists some positive constant $C$ such that
$$u(x)\leq C (1+d^{-\beta_1}(x)),\,\,\,|\nabla u(x)|\leq C (1+d^{-\beta_1-1}(x)),$$
$$v(x)\leq C (1+d^{-\beta_2}(x)),\,\,\,|\nabla v(x)|\leq C (1+d^{-\beta_2-1}(x)),$$
where $\beta_1=\frac{\alpha_2 q_1 +\alpha_1(1-q_2)}{q_1p_2-(1-q_2)(1-p_1)}>0$, $\beta_2=\frac{\alpha_1 p_2 +\alpha_2(1-p_1)}{q_1p_2-(1-q_2)(1-p_1)}>0$.
\label{lem3}
\end{mlem}
Further, we prove
\begin{mthm} \label{th2}
 Suppose $(A_1)-(A_3)$ hold and $0<\theta< \underset{i=1,2}{\min} \{\beta_i, \alpha_i-1\}$ and $\theta(p_i+q_i-1)\leq 1$.
For $1<\alpha_1,\,\alpha_2\,<2,$ if $u\,\in L_{\alpha_1}\cap C^{1,1}_{loc}(\Omega)$ and $v\,\in L_{\alpha_2}\cap C^{1,1}_{loc}(\Omega)$
are positive solutions of system (\ref{In-1.2}), then  for some positive constant $c$ it holds
$$\|u\|_{1,\theta}\leq c, \quad \|v\|_{1,\theta}\leq c.$$
\end{mthm}

As a particular case of  (\ref{In-1.2}), for
\begin{equation}\left\{\begin{array}{ll}
(- \Delta)^{\frac{\alpha_1}{2}}u(x)=f(x,u,v),& \mbox {in\,\,\,} \Omega,\\
(- \Delta)^{\frac{\alpha_2}{2}}v(x)=g(x,u,v),& \mbox {in\,\,\,}
\Omega,\\
u=v=0 & \mbox {in\,\,\,} \mathbb{R}^n\setminus\Omega,
\end{array} \right.\label{In-1.1}
 \end{equation}
we obtain both the a priori estimate of the positive solution and its existence.
\begin{mthm} \label{th1}
 For $0<\alpha_1,\,\alpha_2<1$, suppose $(A_2)-(A_3)$ hold, and $u\in
L_{\alpha_1}\cap C^{1,1}_{loc}(\Omega)$ and $v\in L_{\alpha_2}\cap
C^{1,1}_{loc}(\Omega)$ are upper semi-continuous on
$\bar{\Omega}$. If $(u,\,v)$ is a pair of positive solutions for
(\ref{In-1.1}), then there exists a positive constant $C $ such that $\|u\|_{L^\infty(\Omega)}\leq C$, $\|v\|_{L^\infty(\Omega)}\leq C$.
\end{mthm}

\begin{mthm} \label{20205254}
For $0<\alpha_1,\,\alpha_2<1$, suppose $(A_2)-(A_4)$ hold. Then
(\ref{In-1.1}) has a pair of solutions $(u, \, v) \in
\big(L_{\alpha_1}\cap C^{1,1}_{loc}(\Omega) \big)\times
\big( L_{\alpha_2}\cap
C^{1,1}_{loc}(\Omega)\big)$.
\end{mthm}

Throughout the paper, we use $c$, $C$ and $C_i$, $i \in N$ to denote positive constants whose values may vary from place to place.

\section{Preliminary Results}

Before we show the proofs of the main results, for the convenience of the readers, we list some useful propositions here.

\begin{pro}(Theorem 2.5 in \cite{QX}) Let $g$ be bounded in $R^n\setminus \Omega$ and $f\in C^{\gamma}_{loc}(\Omega)$. Suppose $u$
is a viscosity solution of
$$(-\Delta)^{\frac{\alpha}{2}}u=f\,\,\, \mbox{in }\Omega,\,\,\,u=g\,\,\,\mbox{in }R^n\setminus \Omega.$$
Then  $u\in C^{\alpha+\gamma}_{loc}(\Omega)$.
\label{lem1}
\end{pro}

\begin{pro}\label{lem4}
(Theorem 1.2 in \cite{K}) Assume that $\alpha\in (1,\,2)$. Suppose $u$ is a viscosity solution of
$$(-\Delta)^{\frac{\alpha}{2}}u=f,\,\,\,\mbox{in }\Omega,$$
where $f\in L^{\infty}_{loc}(\Omega)$. Then there exists $\gamma=\gamma(n,\,\alpha)\in (0,\,1)$ such that $u\in C^{1,\,\gamma}_{loc}(\Omega)$.
Moreover, for every ball $B_R\subset \subset \Omega$, there exists a positive constant $C=C(n,\,\alpha,\,R)$ such that
$$\|u\|_{C^{1,\gamma}(\bar{B}_{R/2})} \leq C (\|f\|_{L^{\infty}(B_R)}+\|u\|_{L^{\infty}(R^n)}).$$
\end{pro}

\begin{pro}
(Lemma 5 in \cite{CS}) Let $\{u_k\}$, $k\in {N}$ be a sequence of functions that are bounded in $R^n$ and continuous in
$\Omega$, $f_k$ and $f$ are continuous in $\Omega$ such that

(1) $\Delta^{\alpha}u_k\leq f_k$ in $\Omega$ in viscosity sense.

(2) $u_k\rightarrow u$ locally in $\Omega$.

(3) $u_k\rightarrow u$ a.e. in $R^n$.

(4) $f_k\rightarrow f$ locally uniformly in $\Omega$.

Then $\Delta^{\alpha}u\leq f$ in $\Omega$ in viscosity sense.
\label{lem2}
\end{pro}

\begin{pro}
(Lemma 6 in \cite{BPMQ}) Suppose $0<\alpha<2$. For every $\tau\in(\frac{\alpha}{2},\,\alpha)$, if $u$ satisfies
$$(-\Delta)^{\frac{\alpha}{2}}u \leq C_1 d^{-\tau}_{\lambda},\,\,\,\mbox{in }\Omega_{\lambda}$$
for some $C_1>0$ with $u=0$ in $R^n\setminus \Omega$, then
$$u(x)\leq C_2(C_1+\|u\|_{L^{\infty}(\Omega_{\lambda})})d^{\alpha-\tau}_{\lambda}\,\,\,\mbox{for }x\in (\Omega_{\lambda})_{\delta}$$
for some $C_2>0$ depending on $\alpha$, $\delta$, $\tau$.
\label{lem5}
\end{pro}

\begin{pro}
(Lemma 5 in \cite{BPMQ})
Suppose $1<\alpha<2$. Let $f \in C(\Omega)$ be such that $\|f\|_{0,\,\tau} < +\infty$ for some $\tau \in (\alpha/2,\,\alpha)$. Then the unique
solution of
\begin{equation}
(-\Delta)^{\frac{\alpha}{2}}u=f\,\,\,\mbox{in }\Omega,\,\,\,u=0 \,\,\,\mbox{in }R^n\setminus \Omega
\end{equation}
verifies
$$\|\nabla u\|_{0,\,\tau-\alpha+1}\leq C_0(\|f\|_{0,\,\tau}+\|u\|_{0,\,\tau-\alpha}),$$
where $C_0$ is a positive constant that depends on $n$ and $\alpha$ but not on $\Omega$.
\label{lem6}
\end{pro}

Applying the results above, we derive our key lemma.

\textbf{Proof of Lemma \ref{lem3}.} Assume on contrary, then there exists a sequence of  solutions $\{(u_k,\,v_k)\}$ of  (\ref{lem-2-1}) such that for
$$w_k(x):=u^{\frac{1}{\beta_1}}_k(x)+|\nabla u_k(x)|^{\frac{1}{\beta_1+1}} +  v^{\frac{1}{\beta_2}}_k(x) + |\nabla v_k(x)|^{\frac{1}{\beta_2+1}},$$
and $\{y^k\}\subseteq  \Omega$ such that
\begin{equation}
w_k(y^k)> 2k (1+d^{-1}(y^k)).
\label{lem-2-2}
\end{equation}

By Lemma 5.1 in \cite{PQS}, there exists $\{x^k\}\subseteq  \Omega$ such that
$$w_k(x^k)\geq w_k(y^k),\,\,\,w_k(x^k)>2kd^{-1}(x),$$
\begin{equation}
w_k(z)\leq 2 w_k(x^k),\,\,\,\forall z\in B(x^k,\,kw^{-1}_k(x^k)),\label{lem-2-3}
\end{equation}
where $B(x^k,\,kw^{-1}_k(x^k))$ represents the ball with radius $kw^{-1}_k(x^k)$ centered at $x^k$.

Obviously,
$$w_k(x^k)\rightarrow + \infty ,\,\,\, \mbox{as }k\rightarrow \infty.$$

Define $$\lambda_k:=w^{-1}_k(x^k), \quad \Omega_k:= \{x \mid  \lambda_k x+x^k \in \Omega \},$$
$$\bar{u}_k(x):=\lambda^{\beta_1}_k u_k(\lambda_kx+x^k),\,\,\,x\in \Omega_k,$$
$$\bar{v}_k(x):=\lambda^{\beta_2}_k v_k(\lambda_kx+x^k),\,\,\,x\in \Omega_k.$$
By an elementary calculation, we get
\begin{equation}\left\{\begin{array}{ll}
(- \Delta)^{\frac{\alpha_1}{2}}\bar{u}_k(x)+ \lambda^{\alpha_1-1}_k \sum\limits^n_{i=1}b_i(\lambda_k x+x^k)\frac{\partial \bar{u}_k(x)}{\partial x_i}+ \lambda^{\alpha_1}_k B(\lambda_k x+x^k)\bar{u}_k(x)\\
=\lambda^{\alpha_1+\beta_1}_k f(\lambda_kx+x^k,\lambda^{-\beta_1}_k \bar{u}_k(x),\lambda^{-\beta_2}_k \bar{v}_k(x)),\,\,\,\qquad \mbox {  in } \Omega_k,\\
(- \Delta)^{\frac{\alpha_2}{2}}\bar{v}_k(x)+ \lambda^{\alpha_2-1}_k \sum\limits^n_{i=1}c_i(\lambda_k x+x^k)\frac{\partial \bar{v}_k(x)}{\partial x_i}+ \lambda^{\alpha_2}_k C(\lambda_k x+x^k)\bar{v}_k(x)\\
=\lambda^{\alpha_2+\beta_2}_k g(\lambda_kx+x^k,\lambda^{-\beta_1}_k \bar{u}_k(x),\lambda^{-\beta_2}_k \bar{v}_k(x)),\,\,\,\qquad\mbox {  in } \Omega_k,\\
\end{array} \right.\label{lem-2-4} \end{equation}

For $k$ large enough, with (\ref{lem-2-3}) we derive
\begin{eqnarray}
&&\bar{u}^{\frac{1}{\beta_1}}_k(x)+|\nabla \bar{u}_k(x)|^{\frac{1}{\beta_1+1}} + \bar{v}_k^{\frac{1}{\beta_2}}(x) +
|\nabla \bar{v}_k(x)|^{\frac{1}{\beta_2+1}} \nonumber\\
&=& \frac{w_k(\lambda_k x+x^k)}{w_k(x^k)}\label{lem-2-5}  \\
&\leq& 2.\nonumber
\end{eqnarray}
Therefore $\bar{u}_k$, $\nabla \bar{u}_k$, $\bar{v}_k$, $\nabla \bar{v}_k$ are uniformly bounded in $B_k$.
From Proposition \ref{lem4} we know there exist $\beta', \, \beta''\in (0,\,1)$  such that
$$\bar{u}\in C^{1,\beta'}_{loc}(\Omega), \quad \bar{v}\in C^{1,\beta''}_{loc}(\Omega).$$ Moreover, for every ball $B_R\subset\subset\Omega$, there exists a positive constant $C=C(n,\,\alpha_1,\,\alpha_2,\,R)$ such that
\begin{equation}
\|\bar{u}\|_{C^{1,\beta'}(\bar{B}_{R/2})}\leq C(\|f\|_{L^{\infty}(B_R)}+\|\bar{u}\|_{L^{\infty}(R^n)})\leq C,
\label{lem-2-7}
\end{equation}
\begin{equation}
\|\bar{v}\|_{C^{1,\beta''}(\bar{B}_{R/2})}\leq C(\|g\|_{L^{\infty}(B_R)}+\|\bar{v}\|_{L^{\infty}(R^n)})\leq C.
\label{lem-2-8}
\end{equation}

Using the  diagonal argument and the Arzel\'{a}-Ascoli theorem, we derive that a subsequence of $\{(\bar{u}_k,\,\bar{v}_k)\}$ (still denoted by $\{(\bar{u}_k,\,\bar{v}_k)\}$), satisfies
$$\bar{u}_k\rightarrow \bar{u},\,\,\, \bar{v}_k\rightarrow \bar{v} \mbox{  in } C^1_{loc}(R^n), \quad k\rightarrow \infty.$$
Recall from (\ref{lem-2-5}) we have
\begin{equation}
\bar{u}^{\frac{1}{\beta_1}}_k(0)+|\nabla \bar{u}_k(0)|^{\frac{1}{\beta_1+1}} + \bar{v}_k^{\frac{1}{\beta_2}}(0) +
|\nabla \bar{v}_k(0)|^{\frac{1}{\beta_2+1}}=1.
\label{lem-2-6}
\end{equation}
Taking the limit gives
\begin{equation}
{\bar{u}}^{\frac{1}{\beta_1}}(0)+|\nabla {\bar{u}}(0)|^{\frac{1}{\beta_1+1}} + {\bar{v}}^{\frac{1}{\beta_2}}(0) +
|\nabla {\bar{v}}(0)|^{\frac{1}{\beta_2+1}}=1.
\label{lem-2-9}
\end{equation}
This implies that the limit equation of (\ref{lem-2-4}) possesses nontrivial solution.

Next we derive the nonexistence of positive solution to obtain contradictions with the existence result above.

{\bf{Case i. }}At least one of $\{u_k (x^k)\}$ and $ \{v_k(x^k)\}$ goes to infinity as $k\rightarrow\infty$. Without loss of generality, we may assume
$v_k(x^k)\rightarrow \infty$  as $k \rightarrow \infty$.
Since $\lambda_k\rightarrow 0 $ as $k \rightarrow \infty$, we have
$$v_k( \lambda_k x+x^k) \rightarrow \infty, \quad k \rightarrow \infty.$$
Immediately it yields
$$\lambda^{-\beta_2}_k\bar{v}_k(x)\rightarrow \infty,\,\,\,\mbox{as }k\rightarrow \infty.$$
For $k$ sufficiently large, from condition $(A_3)$, we deduce
\begin{eqnarray}
&& \lambda^{\alpha_1+\beta_1}_k f(\lambda_kx+x^k,\lambda^{-\beta_1}_k \bar{u}_k(x),\lambda^{-\beta_2}_k \bar{v}_k(x)) \nonumber\\
&=& \lambda^{\alpha_1+\beta_1} \frac{f(\lambda_k x+x^k,\,\lambda^{-\beta_1}_{k}\bar{u}_k,\,\lambda^{-\beta_2}_{k}\bar{v}_k)}
{[\lambda^{-\beta_1}_k \bar{u}_k(x)]^{p_1} [\lambda^{-\beta_2}_k \bar{v}_k(x)]^{q_1}}
\, \lambda^{-\beta_1p_1-\beta_2q_1}_k \,\bar{u}^{p_1}_k(x) \, \bar{v}^{q_1}_k(x)\nonumber\\
&=& \lambda^{\alpha_1+\beta_1-\beta_1p_1-\beta_2q_1} \big[k_1(\lambda_k x+ x^k)+ o(1)\big]\, \bar{u}^{p_1}_k(x)\, \bar{v}^{q_1}_k(x)\nonumber\\
&=& k_1(\lambda_kx+x^k) \bar{u}^{p_1}_k(x) \bar{v}^{q_1}_k(x)+o(1).\label{lem-2-10}
\end{eqnarray}
Similarly, we obtain
\begin{eqnarray}
&& \lambda^{\alpha_2+\beta_2}_k g(\lambda_kx+x^k,\lambda^{-\beta_1}_k \bar{u}_k(x),\lambda^{-\beta_2}_k \bar{v}_k(x)) \nonumber\\
&&\qquad\qquad= h_1(\lambda_kx+x^k) \bar{u}^{p_2}_k(x) \bar{v}^{q_2}_k(x)+o(1). \label{lem-2-11}
\end{eqnarray}

Since $\Omega$ is bounded, we know there exists an $x^0 \in \bar{\Omega}$ such that
$$x^k\rightarrow x^0,\,\,\,\mbox{as }k\rightarrow \infty.$$
Then with (\ref{lem-2-10}), (\ref{lem-2-11}) and assumption $(A_1)$,
from Proposition \ref{lem2} we derive the limit equation of (\ref{lem-2-4}) as
\begin{equation}\left\{\begin{array}{ll}
(- \Delta)^{\frac{\alpha_1}{2}}{u}(x)=k_1(x^0) {u}^{p_1}(x) {v}^{q_1}(x),& \mbox {in\,\,\,} \mathbb{R}^n,\\
(- \Delta)^{\frac{\alpha_2}{2}}{v}(x)=h_1(x^0) {u}^{p_2}(x) {v}^{q_2}(x),& \mbox {in\,\,\,} \mathbb{R}^n.
\end{array} \right.\label{lem-2-12} \end{equation}
On the other hand, from the Liouville-type theorem in \cite{ZL} one knows  (\ref{lem-2-12}) has no positive solution. We have a contradiction here. This proves Case i does not exist.

{\bf{Case ii. } }If both $\{ u_k (x^k)\}$ and $\{ v_k(x^k)\}$ are bounded for any $k \in N$, then for some constant $C>0$, it holds
$$u_k(\lambda_kx+x^k)\leq C,\quad \,v_k(\lambda_kx+x^k)\leq C, \quad k \in N.$$
Then
\begin{equation}
\lambda^{-\beta_1}_k\bar{u}_k(x)=u_k(\lambda_kx+x^k)\leq C,\,\,\,\lambda^{-\beta_2}_k\bar{v}_k(x)=v_k(\lambda_kx+x^k)\leq C.
\label{lem-sup-1}
\end{equation}
Combining this with condition $(A_2)$, we arrive at
\begin{equation*}
\lambda^{\alpha_1+\beta_1}_k f(\lambda_kx+x^k,\lambda^{-\beta_1}_k \bar{u}_k(x),\lambda^{-\beta_2}_k \bar{v}_k(x)) \rightarrow 0,\,\,\,\mbox{as }k\rightarrow \infty, \label{lem-2-13}
\end{equation*}
\begin{equation*}
\lambda^{\alpha_2+\beta_2}_k g(\lambda_kx+x^k,\lambda^{-\beta_1}_k \bar{u}_k(x),\lambda^{-\beta_2}_k \bar{v}_k(x)) \rightarrow 0,\,\,\,\mbox{as }k\rightarrow \infty. \label{lem-2-14}
\end{equation*}

Through a similar argument for (\ref{lem-2-12}), one can show that $(u,\,v)$ satisfies
\begin{equation}\left\{\begin{array}{ll}
(- \Delta)^{\frac{\alpha_1}{2}}{u}(x)=0,& \mbox {in\,\,\,} \mathbb{R}^n,\\
(- \Delta)^{\frac{\alpha_2}{2}}{v}(x)=0,& \mbox {in\,\,\,} \mathbb{R}^n.
\end{array} \right.\label{lem-2-15} \end{equation}

Applying the Theorem 2 in \cite{ZCCY} (also see \cite{BKN}) to (\ref{lem-2-15}), we have
$$u\equiv C_3\geq 0,\,\,\,v\equiv C_4\geq 0,\,\,\,\mbox{in }\mathbb{R}^n.$$
Meanwhile it follows from (\ref{lem-sup-1}) that
$$\bar{u}_k \rightarrow 0,\,\,\,\bar{v}_k \rightarrow 0,\,\,\,\mbox{as }k\rightarrow \infty.$$
This implies that
$$u\equiv0,\,\,\,v\equiv0.$$

Hence we obtain a contradiction. This completes the proof of Lemma \ref{lem3}.

\section{The A Priori Estimate}

{\bf{Proof of Theorem \ref{th2}.}}
Suppose otherwise, then there exists a sequence of positive solutions $\{(u_k,\,v_k)\}$ of (\ref{In-1.2}) such that
$$\max\{\|u_k\|_{1,\theta},\,\|v_k\|_{1,\theta}\}\rightarrow \infty,\,\,\,\mbox{as }k\rightarrow \infty.$$
Without loss of generality, we assume
\begin{equation}
(\|u_k\|_{1,\theta})^{\beta_2-\theta}\geq (\|v_k\|_{1,\theta})^{\beta_1-\theta},
\label{th2-1}
\end{equation}
with $\beta_1$ and $\beta_2$ from Lemma {\ref{lem3}}.
Set
$$M_k(x):=d^{\theta}(x)u_k(x)+d^{\theta+1}(x)|\nabla u_k(x)|,$$
$$N_k(x):=d^{\theta}(x)v_k(x)+d^{\theta+1}(x)|\nabla v_k(x)|.$$
Then there exists a sequence $\{x^k\}\subset \Omega$ and $x^0\in \bar{\Omega}$ such that
$$M_k(x^k)\rightarrow +\infty,\,\,\,x^k\rightarrow x^0,\,\,\,\mbox{as }k\rightarrow \infty.$$

Let $\xi^k$ be the projection of $x^k$ on $\partial\Omega$ and set
$$\lambda_k:=(\|u_k\|_{1,\theta})^{-\frac{1}{\beta_1-\theta}}.$$
 Denote
$$\bar{u}_k(x):=\lambda^{\beta_1}_ku_k(\lambda_k x+\xi^k),\quad
\bar{v}_k(x):=\lambda^{\beta_2}_kv_k(\lambda_k x+\xi^k),$$
$$\Omega_k:=\{x\in R^n\,|\,\lambda_kx+\xi^k\in\Omega\}.$$

From (\ref{In-1.2}) we  obtain
\begin{equation}
\left\{\begin{array}{ll}
(- \Delta)^{\frac{\alpha_1}{2}}\bar{u}_k(x)+ \lambda^{\alpha_1-1}_k \sum\limits^n_{i=1}b_i(\lambda_k x+\xi^k)\frac{\partial \bar{u}_k(x)}{\partial x_i}+ \lambda^{\alpha_1}_k B(\lambda_k x+\xi^k)\bar{u}_k(x)\\
\quad\quad\quad \qquad=\lambda^{\alpha_1+\beta_1}_k f\big(\lambda_kx+\xi^k,\lambda^{-\beta_1}_k \bar{u}_k(x),\lambda^{-\beta_2}_k \bar{v}_k(x)\big), \\
(- \Delta)^{\frac{\alpha_2}{2}}\bar{v}_k(x)+ \lambda^{\alpha_2-1}_k \sum\limits^n_{i=1}c_i(\lambda_k x+\xi^k)\frac{\partial \bar{v}_k(x)}{\partial x_i}+ \lambda^{\alpha_2}_k C(\lambda_k x+\xi^k)\bar{v}_k(x)\\
\qquad \quad\qquad=\lambda^{\alpha_2+\beta_2}_k g\big(\lambda_kx+\xi^k,\lambda^{-\beta_1}_k \bar{u}_k(x),\lambda^{-\beta_2}_k \bar{v}_k(x)\big), \\
\qquad  \quad \bar{u}_k(x), \, \bar{v}_k(x)>0,\qquad \quad\quad \qquad \qquad\qquad\qquad x\in \Omega_k,\\
\qquad \quad \bar{u}_k(x)=\bar{v}_k(x)=0, \qquad \qquad\qquad \quad\qquad\quad x \in \mathbb{R}^n\setminus \Omega_k.
\end{array} \right.\label{th2-2} \end{equation}
Let
$$d_k(x):=dist(x,\,\partial \Omega_k).$$
By an elementary calculation, for $x\in \Omega_k$ we derive
\begin{eqnarray}
&&d^{\theta}_k(x)\bar{u}_k(x)+d^{\theta+1}_k(x)|\nabla\bar{u}_k(x)| \nonumber\\
&=& \lambda^{\beta_1-\theta}_k M_k(\lambda_k x+\xi^k) \nonumber\\
&=& \frac{M_k(\lambda_k x+\xi^k)}{\|u_k\|_{1,\theta}}\leq 1.\label{th2-3}
\end{eqnarray}
Similarly, we have
\begin{equation}
d^{\theta}_k(x)\bar{v}_k(x)+d^{\theta+1}_k(x)|\nabla\bar{v}_k(x)|\leq 1, \quad x\in\Omega_k.
\label{th2-4}
\end{equation}
Choosing $y^k=\frac{x^k-\xi^k}{\lambda_k}$, it gives
\begin{equation}
d^{\theta}_k(y^k)\bar{u}_k(y^k)+d^{\theta+1}_k(y^k)|\nabla\bar{u}_k(y^k)|\rightarrow 1,\,\,\,\mbox{as }k\rightarrow\infty.
\label{th2-5}
\end{equation}
By Lemma {\ref{lem3}}, we obtain
\begin{eqnarray}
&&M_k(x^k) \nonumber\\
&=&d^{\theta}(x^k)u_k(x^k)+d^{\theta+1}(x^k)|\nabla u_k(x^k)| \nonumber\\
&\leq& C \big[d^{\theta}(x^k)(1+d^{-\beta_1}(x^k))+ d^{\theta+1}(x^k) (1+d^{-\beta_1-1}(x^k))\big] \nonumber\\
&\leq& Cd^{\theta}(x^k)(1+d^{-\beta_1}(x^k)). \label{th2-6}
\end{eqnarray}
It then follows that
\begin{eqnarray}
\frac{d(x^k)}{\lambda_k}&=&d(x^k)\|u_k\|_{1, \theta}^{\frac{1}{\beta_1-\theta}} \nonumber\\
&\leq& C d^{1+\frac{\theta}{\beta_1-\theta}}(x^k)(1+d^{-\beta_1}(x^k))^{\frac{1}{\beta_1-\theta}} \nonumber\\
&\leq& C d^{\frac{\beta_1}{\beta_1-\theta}}(x^k)(1+d^{-\frac{\beta_1}{\beta_1-\theta}}(x^k)) \nonumber\\
&\leq& C (1+d^{\frac{\beta_1}{\beta_1-\theta}}(x^k)).  \label{th2-7}
\end{eqnarray}

Next we show that, on a subsequence of $\{x^k\}$, still denoted by
$\{x^k\}$, it holds
\begin{equation}
\lim\limits_{k\rightarrow \infty}\frac{d(x^k)}{\lambda_k}:=d>0.
\label{th2-8}
\end{equation}

For $k$  large enough, through a similar argument as in (\ref{lem-2-10}) we derive from   (\ref{th2-3}), (\ref{th2-4}) and assumption $(A_3)$  that
\begin{eqnarray}
&&\lambda^{\alpha_1+\beta_1}_k f(\lambda_kx+\xi^k,\lambda^{-\beta_1}_k \bar{u}_k(x),\lambda^{-\beta_2}_k \bar{v}_k(x)) \nonumber\\
&=& k_1(\lambda_k x +\xi^k) \bar{u}^{p_1}_k(x) \, \bar{v}^{q_1}_k(x)+o(1)
 \nonumber\\
&\leq & C d^{-\theta(p_1+q_1)}_k(x), \label{th2-9}
\end{eqnarray}
and
\begin{equation}
|\lambda^{\alpha_1-1}_k \sum\limits^n_{i=1}b_i(\lambda_k x+\xi^k)\frac{\partial \bar{u}_k(x)}{\partial x_i}|\leq C\lambda^{\alpha_1-1}_k d^{-\theta-1}_k(x),
\label{th2-10}
\end{equation}
and
\begin{equation}
|\lambda^{\alpha_1}_k B(\lambda_k x+\xi^k)\bar{u}_k(x)|\leq C \lambda^{\alpha_1}_k d^{-\theta}_k(x).
\label{th2-11}
\end{equation}

Combining (\ref{th2-9}) through (\ref{th2-11}), for a fixed  small $\delta>0$,
we derive
\begin{equation}
(- \Delta)^{\frac{\alpha_1}{2}}\bar{u}_k(x) \leq C d^{-\theta-1}_k(x),\,\,\,\mbox{ when }d_k(x)<\delta.
\label{th2-12}
\end{equation}

Employing Proposition {\ref{lem5}},
we obtain
\begin{equation}
\bar{u}_k(x)\leq C d^{-\theta-1+\alpha_1}_k(x),\,\,\,\mbox{ when }d_k(x)<\delta.
\label{th2-13}
\end{equation}

By (\ref{th2-3}), we have
\begin{equation}
d^{\theta}_k(x)\bar{u}(x)\leq d^{1-\alpha_1}_k(x) d^{\alpha_1-1}_k(x)\leq \delta^{1-\alpha_1} d^{\alpha_1-1}_k(x),\,\,\,\mbox{ when }d_k(x)\geq\delta.
\end{equation}
That is,
\begin{equation}
\bar{u}_k(x)\leq \delta^{1-\alpha_1} d^{-\theta-1+\alpha_1}_k(x),\,\,\,\mbox{ when }d_k(x)\geq\delta.
\label{th2-14}
\end{equation}
From (\ref{th2-13}) and (\ref{th2-14}), we deduce that
\begin{equation}
\bar{u}_k(x)\leq C d^{-\theta-1+\alpha_1}_k(x).
\label{th2-15}
\end{equation}
Then it follows from Proposition \ref{lem6}  that
\begin{equation}
|\nabla\bar{u}_k(x)|\leq C d^{-\theta-2+\alpha_1}_k(x).
\label{th2-16}
\end{equation}

Combining (\ref{th2-5}), (\ref{th2-15}) and (\ref{th2-16}), at $y_k$ it holds \begin{equation}
d^{\alpha_1-1}_k(y_k)\geq C>0.
\label{2020611}
\end{equation}
It instantly follows from $\alpha_1>1$ that
$$\frac{d(x^k)}{\lambda_k}\geq C>0.$$ Now it is easy to see that
$d>0$ in (\ref{th2-8}).

Meanwhile, from (\ref{2020611}) we can see that $y^k$ is away from $\Omega_k$. Thus there exists some $y^0$ and a subsequence, still denoted by $\{y^k\}$, such that
$$y^k\rightarrow y^0 \in \Omega_k, \quad k \rightarrow \infty. $$
 Using Proposition \ref{lem4}, Arzel\'{a}-Ascoli theorem and a diagonal argument, we can verify
there exist $\bar{u}, \,\bar{v} $ such that
\begin{equation}
\bar{u}_k\rightarrow \bar{u},\,\,\,\bar{v}_k\rightarrow \bar{v},\,\,\,\mbox{in }C^1_{loc}(R^n_{+,\,d}),
\end{equation}
where $R^n_{+,\,d}:=\{x\in R^n\,|\,x_n\geq -d\,\}$ with the $d>0$ from (\ref{th2-8}).

It now follows from (\ref{th2-5}) that
$$dist^{\theta}(y^0,\,R^n_{+,\,d})\bar{u}(y^0)+ dist^{\theta+1}(y^0,\,R^n_{+,\,d})|\nabla\bar{u}(y^0)|=1.$$
This implies that $$(\bar{u},\,\bar{v})\neq (0,\,0).$$

When $dist(x,\,R^n_{+,\,d})<\delta$, from
$$\bar{u}(x)\leq C dist^{-\theta-1+\alpha_1}(x,\,R^n_{+,\,d}),\,\,\,\bar{v}(x)\leq C dist^{-\theta-1+\alpha_1}(x,\,R^n_{+,\,d}),$$
 we know $\bar{u},\,\bar{v}\in C(R^n)$.

Employing Proposition \ref{lem2},  we take limit in (\ref{th2-2}) and arrive at
\begin{equation}\left\{\begin{array}{ll}
(- \Delta)^{\frac{\alpha_1}{2}}\bar{u}(x)=k_1(x^0) \bar{u}^{p_1}(x) \bar{v}^{q_1}(x),& \mbox {in\,\,\,} \mathbb{R}^n_{+,d},\\
(- \Delta)^{\frac{\alpha_2}{2}}\bar{v}(x)=h_1(x^0) \bar{u}^{p_2}(x) \bar{v}^{q_2}(x),& \mbox {in\,\,\,} \mathbb{R}^n_{+,d},\\
\bar{u}=\bar{v}=0 & \mbox {in\,\,\,}\mathbb{R}^n\setminus \mathbb{R}^n_{+,d}.
\end{array} \right.\label{th2-17} \end{equation}
However, by Theorem 2 in \cite{ZL1} we know (\ref{th2-17}) has no positive solution. A contradiction. This completes  the proof.

In the following we show the $L^{\infty}$-norm of solutions for system (\ref{In-1.1}) are bounded.

{\bf{Proof of Theorem \ref{th1}.}} Suppose in the contrary, there exists a sequence of positive solutions  $\{(u_k,\,v_k)\}$ to system (\ref{In-1.1}) such that
$$\max\{\|u_k\|_{L^{\infty}(\Omega)},\,\|v_k\|_{L^{\infty}(\Omega)}\}\,\rightarrow\,\infty,\,\,\,\mbox{as }k\rightarrow \infty.$$

Without any loss of generality, we assume
$$\|u_k\|^{\beta_2}_{L^{\infty}(\Omega)}\,\geq\,\|v_k\|^{\beta_1}_{L^{\infty}(\Omega)}.$$
Denote
$$\lambda_k\,:=\,\|u_k\|^{-\frac{1}{\beta_1}}_{L^{\infty}(\Omega)}.$$
It is easy to see that
$$\lambda_k\,\rightarrow\,0,\,\,\,\mbox{as }k\,\rightarrow\,\infty.$$
Let $x^k\,\in\,\Omega$ be a point where $u_k$ assumes its maximum. Define
$$\widetilde{u}_k(x):=\lambda^{\beta_1}_k u_k(\lambda_kx+x^k),\,\,\,\widetilde{v}_k(x):=\lambda^{\beta_2}_k v_k(\lambda_kx+x^k).$$
Then we have
$$\widetilde{u}_k(0)=1, \quad 0\leq \widetilde{u}_k,\,\widetilde{v}_k\leq c, \quad x \in\Omega.$$ By an elementary calculation, we derive that
$(\widetilde{u}_k,\,\widetilde{v}_k)$ satisfies the following system
\begin{equation}\left\{\begin{array}{ll}
(- \Delta)^{\frac{\alpha_1}{2}}\widetilde{u}_k(x)=
\lambda^{\alpha_1+\beta_1}_{k}f(\lambda_k x+x^k,\,\lambda^{-\beta_1}_{k}\widetilde{u}_k,\,\lambda^{-\beta_2}_{k}\widetilde{v}_k),& \mbox {in\,\,\,} \Omega_k,\\
(- \Delta)^{\frac{\alpha_2}{2}}\widetilde{v}_k(x)=
\lambda^{\alpha_2+\beta_2}_{k}g(\lambda_k x+x^k,\,\lambda^{-\beta_1}_{k}\widetilde{u}_k,\,\lambda^{-\beta_2}_{k}\widetilde{v}_k),& \mbox {in\,\,\,}
\Omega_k,\\
\widetilde{u}_k>0,\,\,\,\widetilde{v}_k>0,& \mbox {in\,\,\,}
\Omega_k,\\
\widetilde{u}_k=\widetilde{v}_k=0 & \mbox {in\,\,\,} \mathbb{R}^n\setminus\Omega_k,
\end{array} \right.\label{est-1.3} \end{equation}
where $\Omega_k:=\{x\in \mathbb{R}^n\,|\,\lambda_k x+ x^k \in \Omega\}.$

By assumption $(A_3)$, similar to (\ref{th2-9}) we have
\begin{eqnarray}
&&\lambda^{\alpha_1+\beta_1}f(\lambda_k x+x^k,\,\lambda^{-\beta_1}_{k}\widetilde{u}_k,\,\lambda^{-\beta_2}_{k}\widetilde{v}_k)\nonumber\\
&=& k_1(\lambda_k x +x^k) \widetilde{u}^{p_1}_k(x) \, \widetilde{v}^{q_1}_k(x)+o(1),\label{est-1.4}
\end{eqnarray}
and
\begin{eqnarray}
&&\lambda^{\alpha_2+\beta_2}g(\lambda_k x+x^k,\,\lambda^{-\beta_1}_{k}\widetilde{u}_k,\,\lambda^{-\beta_2}_{k}\widetilde{v}_k)\nonumber\\
&=& h_1(\lambda_k x +x^k) \widetilde{u}^{p_2}_k(x) \, \widetilde{v}^{q_2}_k(x)+o(1).\label{est-1.5}
\end{eqnarray}

Let $d_k(x):=dist(x^k,\,\partial\Omega)$. We use contradiction argument to
prove our result.

{\bf{Case i) }} $\lim\limits_{k\rightarrow \infty}\frac{d_k}{\lambda_k}=\infty.$

It is easy to see that
$$\Omega_k \rightarrow \mathbb{R}^n,\,\,\,\mbox{as }k\rightarrow \infty.$$

Since $\widetilde{u}_k(x)$ and $\widetilde{v}_k(x)$ are uniformly bounded, by Proposition \ref{lem1} and the Arzel\'{a}-Ascoli theorem,  we derive
\begin{equation}
\widetilde{u}_k\rightarrow \widetilde{u},\,\,\,\widetilde{v}_k\rightarrow \widetilde{v},\,\,\,\mbox{in }C^1_{loc}(\mathbb{R}^n).
\label{est-1.6}
\end{equation}

There exists a subsequence of $\{x^k\}$, also denoted by $\{x^k\}$, such that $x^k\rightarrow x^0$ in $\bar{\Omega}$.  Proposition \ref{lem2} implies that $(\widetilde{u},\,\widetilde{v})$ is a positive solution of
\begin{equation}\left\{\begin{array}{ll}
(- \Delta)^{\frac{\alpha_1}{2}}\widetilde{u}(x)=k_1(x^0) \widetilde{u}^{p_1}(x) \widetilde{v}^{q_1}(x),& \mbox {in\,\,\,} \mathbb{R}^n,\\
(- \Delta)^{\frac{\alpha_2}{2}}\widetilde{v}(x)=h_1(x^0) \widetilde{u}^{p_2}(x) \widetilde{v}^{q_2}(x),& \mbox {in\,\,\,} \mathbb{R}^n,\\
\end{array} \right.\label{est-1.7} \end{equation}
in the viscosity sense. However, by Theorem 1.1 in \cite{ZL}, we know that (\ref{est-1.7}) has no positive solution when
 $p_1+\frac{n-\alpha_2}{n-\alpha_1}q_1<\frac{n+\alpha_1}{n-\alpha_1}$ and
$ p_2+\frac{ n-\alpha_2}{n-\alpha_1 } q_2< \frac{n+\alpha_2}{n-\alpha_1 }$.
The contradiction proves that Case i does not exist.

{\bf{Case ii) }} $\lim\limits_{k\rightarrow \infty}\frac{d_k}{\lambda_k}=d>0.$

It is not difficult to see that
$$\Omega_k \rightarrow \mathbb{R}^n_{+,d},\,\,\,\mbox{as } k\rightarrow \infty,$$
where $\mathbb{R}^n_{+,d}:=\{x\in R^n\,|\,x_n\geq -d\,\}.$

Similar to Case i, we can prove there exists $(\widetilde{u},\,\widetilde{v})$ such that as
$k\rightarrow \infty$,
\begin{equation}
\widetilde{u}_k\rightarrow \widetilde{u},\,\,\,\widetilde{v}_k\rightarrow \widetilde{v},
\label{est-1.8}
\end{equation}
and
\begin{equation}
(-\Delta)^{\frac{\alpha_1}{2}}\widetilde{u}_k\rightarrow (-\Delta)^{\frac{\alpha_1}{2}}\widetilde{u},
\,\,\,(-\Delta)^{\frac{\alpha_2}{2}}\widetilde{v}_k\rightarrow (-\Delta)^{\frac{\alpha_2}{2}}\widetilde{v}.
\end{equation}

Moreover, $\{x^k\}$ has a subsequence that converges to $ x^0\,\in \partial \Omega$. As a result, we have
\begin{equation}\left\{\begin{array}{ll}
(- \Delta)^{\frac{\alpha_1}{2}}\widetilde{u}(x)=k_1(x^0) \widetilde{u}^{p_1}(x) \widetilde{v}^{q_1}(x),& \mbox {in\,\,\,} \mathbb{R}^n_{+,d},\\
(- \Delta)^{\frac{\alpha_2}{2}}\widetilde{v}(x)=h_1(x^0) \widetilde{u}^{p_2}(x) \widetilde{v}^{q_2}(x),& \mbox {in\,\,\,} \mathbb{R}^n_{+,d},\\
\widetilde{u}=\widetilde{v}=0, & \mbox {in\,\,\,}\mathbb{R}^n\setminus \mathbb{R}^n_{+,d},
\end{array} \right.\label{est-1.9} \end{equation}

In \cite{ZL1} the authors proved that (\ref{est-1.9}) has no positive solution. Meanwhile, by the definition of $\widetilde{u}_k$ and
$\widetilde{v}_k$, we have
$$\widetilde{u}(0)=\lim\limits_{k\rightarrow \infty}\widetilde{u}_k(0)=1.$$
This is a contradiction. Case ii also does not exist.

{\bf{Case iii) }}  $\lim\limits_{k\rightarrow \infty}\frac{d_k}{\lambda_k}=0.$

We may assume a subsequence of $\{x^k\}$, still denoted by $\{x^k\}$, converges to $x^0\,\in \partial\Omega$.

Let $\widetilde{x}^k:=\frac{x^0-x^k}{\lambda_k}$, we have
$$\widetilde{u}_k(\widetilde{x}^k)=\widetilde{v}_k(\widetilde{x}^k)=0,$$
and
$$\widetilde{x}^k\rightarrow 0,\,\,\,\mbox{as }k\rightarrow \infty.$$

Similar to the argument of (2.18) in \cite{CLL}, we can prove
for $x \in \Omega_k$  and close to  $\partial\Omega_k$,
\begin{equation}
|\widetilde{u}_k(x)-\widetilde{u}_k(\widetilde{x}^k)|\leq c |x-\widetilde{x}^k|^{\frac{\alpha_1}{2}}.
\label{est-1.10}
\end{equation}
Thus the value of $ |\widetilde{u}_k(x)-\widetilde{u}_k(\widetilde{x}^k)|$ can be arbitrarily small for $k$ sufficiently large. However,
\begin{equation}
\widetilde{u}_k(0)-\widetilde{u}_k(\widetilde{x}^k)=1.
\label{est-1.11}
\end{equation}
This is a contradiction. Therefore Case iii is impossible.

\section{Existence of Solution}

In this section we use the a priori bounds obtained in Theorem \ref{th1} and the topological degree theorem to prove Theorem \ref{20205254}.

\textbf{Proof of Theorem \ref{20205254}.} Let
 \begin{equation}\label{20205181}
 T(u,v):=( T_{\alpha_1}(u,v), T_{\alpha_2}(u,v))
 \end{equation}
 with
 \begin{equation*}
   T_{\alpha_1}(u,v):= \int_\Omega G_{\alpha_1}(x,y) f(y, u(y), v(y))\,dy,
 \end{equation*}
  \begin{equation*}
   T_{\alpha_2}(u,v):= \int_\Omega G_{\alpha_2}(x,y) g(y, u(y), v(y))\,dy.
 \end{equation*}

\begin{lem}
 The operator $T$ defined in (\ref{20205181}) is compact.
\end{lem}

\emph{Proof.} We first prove the continuity of $T$.

  Take a sequence of functions $\{u_n\}$,$ \{v_n\}$ that converges respectively to $u$ and $v$ in $C^0(\Omega)$.
  For a given $\epsilon \rightarrow 0$, there exists an $N>0$ such that for $n>N$,
  \begin{equation*}
    |u_n-u|<\epsilon, \qquad |v_n-v|<\epsilon.
  \end{equation*}
  It follows from the continuity of $f$ that
    \begin{equation*}
    |f(y,u_n(y), v_n(y))-f(y,u_n(y), v(y))|<\epsilon,
  \end{equation*}
  and
    \begin{equation*}
  | f(y,u_n(y), v(y))-f(y,u(y), v(y))|<\epsilon.
  \end{equation*}
  Hence
  \begin{eqnarray*}
&&|T_{\alpha_1}(u_n,v_n)-T_{\alpha_1}(u,v)|\\
&=&|\int_\Omega G_{\alpha_1}(x,y)[ f(y, u_n(y), v_n(y))-f(y, u(y), v(y))] dy|
\nonumber\\
&\leq& \int_\Omega G_{\alpha_1}(x,y)| f(y, u_n(y), v_n(y))-f(y, u_n(y), v(y))| dy\\
&& +\int_\Omega G_{\alpha_1}(x,y)| f(y, u_n(y), v(y))-f(y, u_n(y), v(y))| dy
\nonumber\\
&\leq& 2\epsilon \int_\Omega G_{\alpha_1}(x,y) \, dy\\
&<& \epsilon.
\end{eqnarray*}
The last inequality is a result of the uniform bounded-ness of $\int_\Omega G_{\alpha_1}(x,y)dy$ in $\Omega$ proved in Lemma 3 in \cite{BPMQ} and Proposition 1.1 in \cite{RS}. This proves $T_{\alpha_1}$ is continuous. Through a similar argument one can prove the continuity of $T_{\alpha_2}$ and derive that $T$ is continuous.

Next we show that $T$ is compact.

Let $\{u_n\}$ and $\{v_n\}$  be bounded sequences in $C^0(\Omega)$. Then
\begin{eqnarray}
|T_{\alpha_1}(u_n,v_n)|&=& |\int_\Omega G_{\alpha_1}(x,y)f(y, u_n(y), v_n(y)) \, dy|\nonumber\\
&\leq & \int_\Omega G_{\alpha_1}(x,y)|f(y, u_n(y), v_n(y)) |\, dy\nonumber\\
&\leq & C \, \underset{\Omega}{\max} |f|\nonumber\\
&\leq & C .
\end{eqnarray}
This shows the bounded-ness of $\{T_{\alpha_1}(u_n,v_n)\}$ in $C^0(\Omega)$. The bounded-ness of $\{T_{\alpha_2}(u_n,v_n)\}$ can be verified in a similar way.

Let $(\bar{u}_n,\bar{v}_n):=T_{\alpha_1}(u_n,v_n)$. Then
$$
\left\{\begin{array}{ll}
(-\lap)^{\alpha_1/2}\bar{u}_n= f(x, u_n, v_n), & x \in \Omega,\\
(-\lap)^{\alpha_2/2}\bar{v}_n= g(x, u_n, v_n), & x \in \Omega.
 \end{array}
 \right.
$$
Since $\bar{u}_n$, $\bar{v}_n$ are bounded in $W^{\alpha_1, \infty}(\Omega)$ and $W^{\alpha_2, \infty}(\Omega)$ respectively, by Sobolev embedding(\cite{AF}) we have
$$W^{\alpha_1, \infty}(\Omega)\hookrightarrow\hookrightarrow C^0(\Omega).$$
This guarantees $\bar{u}_n$ and $\bar{v}_n$ each has a converging subsequence in $C^0(\Omega)$.

We thus prove the compactness of $T$.

 Now we show the existence of position solutions through a combination of the a priori estimate obtained in (\ref{th1})
 and the topological degree theory below.

 \begin{pro}\label{20205191}
Suppose that $(X,P)$ is an ordered Banach space, $U\subset P$ is bounded open and contains $0$. Assume that there exists $ \rho>0$ such that $B_\rho(0)\cap P\subset U$ and that $K:\bar{U}\mapsto P$ is compact and satisfies:
\begin{enumerate}
\item For any $x\in P$ with $|x|=\rho$, and $\lambda\in[0,1)$, $x\neq\lambda Kx$;
\item  There exists some $ y\in P\backslash\{0\}$, such that $\ x-Kx\neq ty$ for any $t\geq 0$ and $x\in \partial U$.
\end{enumerate}
Then $K$ possesses a fixed point on $\bar{U}_{\rho}$, where $U_{\rho}=U\backslash B_{\rho}(0)$.
\end{pro}

We verify that the positive solution of Eq. (\ref{In-1.1}) must satisfy the two conditions in Proposition \ref{20205191}.

First, we show that there exits some $\rho>0$ small such that for any $u,\, v \in \partial B_\rho(0)$, it holds
\begin{equation}\label{20205193}
  (u,v)\neq \lambda T(u,v), \quad  \lambda \in [0,1).
\end{equation}
In fact it suffices to show that for $\rho>0$ small and
\begin{equation}\label{20205192}
\| u\|_{C^0(\Omega)}=\| v\|_{C^0(\Omega)}=\rho,
\end{equation}
it holds that
\begin{equation*}
  \| T_{\alpha_1}(u,v)\|_{C^0(\Omega)}< \| u \|_{C^0(\Omega)}, \quad
    \| T_{\alpha_2}(u,v)\|_{C^0(\Omega)}< \| v\|_{C^0(\Omega)}.
\end{equation*}

For any fixed small $\varepsilon_0>0$ and large $M>0$, let
\begin{eqnarray*}
I&:=& \{x \in \Omega \mid \min\{\| u\|_{C^0(I)}, \| v \|_{C^0(I)}\}<\varepsilon_0\} \\
II &:=& \{x \in \Omega \mid \max\{\| u\|_{C^0(II)}, \| v \|_{C^0(II)}\}>M\}.
\end{eqnarray*}
From (\ref{20205192}), $(A_3)$ and $(A_4)$
we have
\begin{eqnarray*}
0&\leq&T_{\alpha_1}(u,v)=\int_\Omega G_{\alpha_1}(x,y) f(y, u(y), v(y))\,dy\\
&=&\int_{I} G_{\alpha_1}(x,y)
\frac{ f(y, u(y), v(y))}{|u|^{\tau_1}|v|^{\eta_1}}\,|u|^{\tau_1}|v|^{\eta_1}\,dy\\
&& + \int_{II} G_{\alpha_1}(x,y)
\frac{ f(y, u(y), v(y))}{|u|^{p_1}|v|^{q_1}}\,|u|^{p_1}|v|^{q_1}\,dy\\
&& +\int_{\Omega \setminus (I \cup II) }  G_{\alpha_1}(x,y)
\frac{ f(y, u(y), v(y))}{|u|^{\min\{p_1,\tau_1 \}}|v|^{\min\{q_1, \eta_1\}}}\,|u|^{\min\{p_1,\tau_1 \}}|v|^{\min\{q_1, \eta_1\}}\,dy\\
&\leq & C\bigg(\| u\|^{\tau_1}_{C^0(I)}\| v \|^{\eta_1}_{C^0(I)}+ \| u\|^{p_1}_{C^0(II)} \| v \|^{q_1}_{C^0(II)}\\
&&+ \| u\|^{\min\{p_1,\tau_1 \}}_{C^0(\Omega \setminus (I \cup II))}\| v \|^{\min\{q_1, \eta_1\} }_{C^0(\Omega \setminus (I \cup II))}\bigg) \\
&\leq& C\| u\|^{\min\{p_1,\tau_1 \}}_{C^0(\Omega )}\| v \|^{\min\{q_1, \eta_1\} }_{C^0(\Omega)}  \\
&<&  \| u\|_{C^0(\Omega)}.
\end{eqnarray*}
Through similar steps one can prove that
\begin{equation*}
  0\leq T_{\alpha_2}(u,v)<\| v \|_{C^0(\Omega)}.
\end{equation*}
This yields (\ref{20205193}).

Next we show that there exists a pair of functions $(\varphi, \phi) \in (P \times P)$ and
$(\varphi, \phi)\neq(0,0)$ such that
\begin{equation}\label{20205201}
  (u,v)-T(u,v)\neq t(\varphi, \phi), \, \forall t\geq 0, \, \forall  (u,v)\in \partial B_R(0) \cap U.
\end{equation}
Let $\varphi,\, \phi$ be the unique viscosity solution (see \cite{FQ}) to
\begin{multicols}{2}
\begin{equation}
\left\{\begin{array}{ll}
(-\lap)^{\alpha_1/2} \phi (x) =1,  & x \in\Omega,\\
\phi (x) >0,  & x \in\Omega,\\
\phi (x)=0,  & x \not \in\Omega.
\end{array}
\right.
\end{equation}

\begin{equation}
\left\{\begin{array}{ll}
(-\lap)^{\alpha_2/2} \varphi (x) =1,  & x \in\Omega,\\
\varphi (x) >0,  & x \in\Omega,\\
\varphi (x)=0,  & x \not \in\Omega.
\end{array}
\right.
\end{equation}
\end{multicols}
To prove (\ref{20205201}), it suffices to prove that for any $t>0$,
\begin{equation}\label{20205202}
\left\{\begin{array}{ll}
(-\lap)^{\alpha_1/2} u (x) =f(x, u,v)+t,  & x \in\Omega,\\
(-\lap)^{\alpha_2/2} v (x) =g(x, u,v)+t,  & x \in\Omega,\\
u(x), \, v(x) >0,  & x \in\Omega,\\
u (x)=v(x)=0,  & x \not \in\Omega,
\end{array}
\right.
\end{equation}
does not have a solution.

Let
\begin{eqnarray*}
\lambda&:= &\inf \{  I \mid I:= \int_{\Omega} [|(-\lap)^{\alpha_1/4}u|^2+|(-\lap)^{\alpha_2/4}v|^2]dx,\\
&& (u,v) \in H^{\alpha_1/2}_0(\Omega)\times H^{\alpha_2/2}_0(\Omega),
\, \int_{\Omega} u(x)v(x)\, dx=1\}.
\end{eqnarray*}
Denote by $(w_{\alpha_1}(x), w_{\alpha_2}(x))$ at which $\inf I$ is attained. Then
\begin{equation*}
  \lambda=\int_{\Omega} [|(-\lap)^{\alpha_1/4}w_{\alpha_1}|^2+|(-\lap)^{\alpha_2/4}w_{\alpha_2}|^2]dx.
\end{equation*}
Let
\begin{equation*}
  \tau_{\alpha_1}:=\int_{\Omega} |(-\lap)^{\alpha_1/4}w_{\alpha_1}|^2 dx, \quad
  \tau_{\alpha_2}:=\int_{\Omega} |(-\lap)^{\alpha_2/4}w_{\alpha_2}|^2 dx.
\end{equation*}
Then $(w_{\alpha_1}(x), w_{\alpha_2}(x))$ solves
\begin{equation}\label{20202222}
\left\{\begin{array}{ll}
(-\lap)^{\alpha_1/2} w_{\alpha_1} (x) =  \tau_{\alpha_1} w_{\alpha_2}(x),  & x \in\Omega,\\
(-\lap)^{\alpha_2/2} w_{\alpha_2} (x) =\tau_{\alpha_2} w_{\alpha_1}(x),  & x \in\Omega,\\
w_{\alpha_1} (x)=w_{\alpha_2}(x)=0,  & x \not \in\Omega.
\end{array}
\right.
\end{equation}

We also have the following result.
\begin{lem}\label{20205211}
\begin{equation*}
w_{\alpha_1}(x),\,w_{\alpha_2}(x)\geq 0,  \quad  x  \in\Omega.
\end{equation*}
\end{lem}

\emph{Proof of Lemma \ref{20205211}.}  Let
\begin{eqnarray*}
  \Omega_{\alpha_1}^+:= \{x \in \Omega \mid w_{\alpha_1} (x)>0   \}, &&\Omega_{\alpha_1}^-:= \{x \in \Omega \mid w_{\alpha_1} (x)<0   \} ,\\
\Omega_{\alpha_2}^+:= \{x \in \Omega \mid w_{\alpha_2} (x)>0   \}, &&\Omega_{\alpha_2}^-:= \{x \in \Omega \mid w_{\alpha_2} (x)<0   \}.
\end{eqnarray*}
We claim that
\begin{equation}\label{20205212}
 \Omega_{\alpha_1}^+=\Omega_{\alpha_2}^+,
\end{equation}
or, equivalently,
\begin{equation}\label{20205221}
  \Omega_{\alpha_1}^+ \cap \Omega_{\alpha_2}^-=  \emptyset.
\end{equation}
If not, then there exists some $\bar{x} \in  \Omega_{\alpha_1}^+ \cap \Omega_{\alpha_2}^-$ such that
\begin{equation*}
  w_{\alpha_1} (\bar{x})=\min _\Omega w_{\alpha_1} (x),\quad \,w_{\alpha_2}(\bar{x})>0.
\end{equation*}
Hence
\begin{eqnarray*}
 (-\lap)^{\alpha_1/2} w_{\alpha_1} (\bar{x})   &=& CP.V. \int_{\mathbb{R} } \frac{w_{\alpha_1} (\bar{x})-w_{\alpha_1} (y)}{|\bar{x}-y|^{n+\alpha_1}}\\
  &=& CP.V. \int_{\Omega_{\alpha_1}^+ } \frac{w_{\alpha_1} (\bar{x})-w_{\alpha_1} (y)}{|\bar{x}-y|^{n+\alpha_1}}+CP.V. \int_{\Omega_{\alpha_1}^- } \frac{w_{\alpha_1} (\bar{x})-w_{\alpha_1} (y)}{|\bar{x}-y|^{n+\alpha_1}}\\
  &&+ CP.V. \int_{\mathbb{R} \setminus \Omega } \frac{w_{\alpha_1} (\bar{x})-w_{\alpha_1} (y)}{|\bar{x}-y|^{n+\alpha_1}} \\
  &<&0.
\end{eqnarray*}
On the other hand,
\begin{equation*}
  (-\lap)^{\alpha_1/2} w_{\alpha_1} (\bar{x})=\tau_{\alpha_1} w_{\alpha_2}(\bar{x})>0.
\end{equation*}
The contradiction proves (\ref{20205221}).

By (\ref{20202222}), (\ref{20205212}) and the Parseval's identity, we have
\begin{eqnarray}\label{20205231}
\nonumber
&& \int_\Omega \tau_{\alpha_1} w_{\alpha_2}^+(x)w_{\alpha_1}^+(x)\, dx\\\nonumber
 &=& \int_{\mathbb{R}^n} \tau_{\alpha_1} w_{\alpha_2}(x) w_{\alpha_1}^+(x)\, dx \\\nonumber
   &=& \int_{\mathbb{R}^n}  (-\lap)^{\alpha_1/2}w_{\alpha_1}(x)  w_{\alpha_1}^+(x)\, dx  \\\nonumber
   &=& \int_{\mathbb{R}^n}  (-\lap)^{\alpha_1/4}w_{\alpha_1}(x) (-\lap)^{\alpha_1/4}  w_{\alpha_1}^+(x)\, dx   \\\nonumber
   &=&  \int_{\mathbb{R}^n}  [(-\lap)^{\alpha_1/4}w_{\alpha_1}^+(x) -(-\lap)^{\alpha_1/4}w_{\alpha_1}^-(x)](-\lap)^{\alpha_1/4}  w_{\alpha_1}^+(x)\, dx \\
  &=& \int_{\mathbb{R}^n} \big\{ [(-\lap)^{\alpha_1/4}w_{\alpha_1}^+(x) ]^2-
  (-\lap)^{\alpha_1/4}w_{\alpha_1}^-(x) (-\lap)^{\alpha_1/4}  w_{\alpha_1}^+(x)\big\}\, dx.
\end{eqnarray}
 Similarly,
 \begin{eqnarray}\label{20205232}
 \nonumber
&& \int_\Omega \tau_{\alpha_1} w_{\alpha_2}^-(x)w_{\alpha_1}^-(x)\, dx\\
 &=&\int_{\mathbb{R}^n} \big\{   (-\lap)^{\alpha_1/4}w_{\alpha_1}^+(x) (-\lap)^{\alpha_1/4}  w_{\alpha_1}^-(x)-
 [(-\lap)^{\alpha_1/4}w_{\alpha_1}^-(x) ]^2\big\}\, dx.
 \end{eqnarray}
 Recall
\begin{eqnarray*}
  \tau_{\alpha_1} &=& \int_\Omega \tau_{\alpha_1} w_{\alpha_2}(x)w_{\alpha_1}(x)\, dx \\
   &=&  \int_{\mathbb{R}^n}  [(-\lap)^{\alpha_1/4}w_{\alpha_1}(x) ]^2 \, dx\\
   &=& \int_{\mathbb{R}^n} \big\{ [(-\lap)^{\alpha_1/4}w_{\alpha_1}^+(x) ]^2-
  2(-\lap)^{\alpha_1/4}w_{\alpha_1}^-(x) (-\lap)^{\alpha_1/4}  w_{\alpha_1}^+(x)\\
  &&+
  [(-\lap)^{\alpha_1/4}w_{\alpha_1}^-(x) ]^2\big\}\, dx,
\end{eqnarray*}
together with (\ref{20205231}) and (\ref{20205232}), we derive
\begin{equation*}
  1= \int_\Omega (w_{\alpha_2}^+(x)w_{\alpha_1}^+(x) - w_{\alpha_2}^-(x)w_{\alpha_1}^-(x))\, dx.
\end{equation*}
Meanwhile, we have
\begin{eqnarray*}
  1  &=&\int_\Omega w_{\alpha_2}(x)w_{\alpha_1}(x)\, dx\\
   &=&\int_\Omega (w_{\alpha_2}^+(x)w_{\alpha_1}^+(x) + w_{\alpha_2}^-(x)w_{\alpha_1}^-(x))\, dx.
\end{eqnarray*}
Now it is easy to see that
\begin{equation*}
  \int_\Omega  w_{\alpha_2}^-(x)w_{\alpha_1}^-(x))\, dx=0.
\end{equation*}
As a result,
\begin{equation*}
  w_{\alpha_2}^-(x)=w_{\alpha_1}^-(x)=0, \quad x \in \Omega.
\end{equation*}
In other words,
\begin{equation*}
w_{\alpha_1}(x), \,  w_{\alpha_2}(x)\geq0, \quad x \in \Omega.
\end{equation*}
This proves the lemma.

For $\bar{\lambda}>\lambda$, suppose (\ref{20205202}) has a pair of solutions $(u, \, v) $ that are continuous in $\Omega$, then there exists some constant
$c_o>0$ such
\begin{eqnarray*}
  f(x, u, v)-\bar{\lambda} v &\geq& -c_o,  \\
 g(x, u, v)-\bar{\lambda} u &\geq& -c_o.
\end{eqnarray*}
Then
\begin{eqnarray}
   \nonumber
        (-\lap)^{\alpha_1/2} u(x) &=& f(x, u,v)+t \\\label{20205241}
         &\geq & \bar{\lambda} v(x) -c_o+t, \quad x \in \Omega.
      \end{eqnarray}
\begin{itemize}
  \item For $t>c_o$, Eq. (\ref{20205241}) becomes
  \begin{equation}\label{20205251}
     (-\lap)^{\alpha_1/2} u(x)\geq  \bar{\lambda} v(x), \quad x \in \Omega.
  \end{equation}
From Eq. (\ref{20202222}),  we derive
  \begin{eqnarray*}
    \int_\Omega  (-\lap)^{\alpha_1/2} u(x) w_{\alpha_1}(x)\, dx &=&   \int_\Omega  u(x) (-\lap)^{\alpha_1/2} w_{\alpha_1}(x) \, dx \\
     &=&  \int_\Omega  u(x)\tau_{\alpha_1} w_{\alpha_2}(x)\, dx.
     \end{eqnarray*}
  Together with (\ref{20205251}), we have
    \begin{equation}\label{20205252}
\tau_{\alpha_1}\int_\Omega  u(x) w_{\alpha_2}(x)\, dx\geq  \bar{\lambda} \int_\Omega  v(x) w_{\alpha_1}(x)\, dx.
 \end{equation}
 In a similar way, one can prove
  \begin{equation}\label{20205253}
\tau_{\alpha_2}\int_\Omega v(x) w_{\alpha_1}(x) \, dx\geq  \bar{\lambda} \int_\Omega u(x) w_{\alpha_2}(x) \, dx.
 \end{equation}
From (\ref{20205252}) and (\ref{20205253}) it yields
\begin{equation*}
\tau_{\alpha_1 } -\frac{\bar{\lambda}^2}{\tau_{\alpha_2}}\geq 0.
\end{equation*}
This implies that at least one of $\tau_{\alpha_1 },\, \tau_{\alpha_2}$ is bigger than $\bar{\lambda}$, which
contradicts
\begin{equation*}
 \bar{\lambda}  > \lambda=\tau_{\alpha_1} + \tau_{\alpha_2}.
\end{equation*}
Thus we conclude that for $t>c_o$, Eq. (\ref{20205202}) does not possess a solution.

  \item For $t\leq c_o$, since $f,\, g$ are continuous on bounded domain $\Omega$, we know the
  terms on the right-hand side of the PDEs from Eq.(\ref{20205202}) are bounded. Then it
  follows from  Theorem  \ref{th1} that
   the solutions to (\ref{20205202}) must be bounded. In other words, there
  exists some $M_o>0$ such that
  \begin{equation*}
    \| u\|_{C^0(\Omega)},\,  \| v\|_{C^0(\Omega)}\leq M_o.
  \end{equation*}
  Therefore, for $R>M_o $, Eq. (\ref{20205202}) has no solution on $\partial B_R(0) \cap U$.
\end{itemize}
By now we have verified the two conditions in Proposition \ref{20205191}. It then follows
from the proposition that there exists a pair of functions $(u,v) \in U \setminus B_\rho (0)$ such
that
\begin{equation*}
(u,v)=T(u,v).
\end{equation*}
This leads to the existence of solution to Eq. (\ref{In-1.1}). 

\section{Acknowledgement}

The authors would like express their gratefulness to Prof. Wenxiong Chen for offering insightful discussions about this subject.

{\em Authors' Addresses and E-mails:}
\medskip

Ran Zhuo

Department of Mathematical and Statistics

Huanghuai University,

Zhumadian, Henan 463000, P. R. China

zhuoran1986@126.com

\medskip

Yan Li

Department of Mathematics

Baylor University

Waco, TEXAS 76798, USA

Yan\_Li1@baylor.edu

\end{document}